\documentclass[10pt]{article}

\def\today{17/10/08} %PLEASE PUT HERE THE DATE 
\usepackage{amsmath,amsfonts,amssymb}
\usepackage{amsthm,amscd}
\usepackage[dvips]{graphicx}
\usepackage{color}
\renewcommand{\Re}{\mathop{\rm Re}\nolimits}
\renewcommand{\Im}{\mathop{\rm Im}\nolimits}

\def\im{{\rm i}}
\def\norma#1{\left\| #1\right\|}

\def\mod#1{#1}

\def\poisson#1#2{\left\{#1,#2\right\}}

\def\sleq{\leq\kern-6pt \cdot\null\hskip4pt}

\def\uno{\mathbb{I}}
\def\jumpsec{\vskip 75 pt}

\def\smjump{\vskip 25 pt}
\def\opnorma#1{\left\bracevert #1 \right\bracevert\null\kern-5pt}
\def\Ham#1#2{H_#2}

\def\lie{\mathcal{L}}
\def\poly{{\mathcal H}}
\def\tond#1{\left(#1\right)}
\def\quadr#1{\left[#1\right]}

\def\eps{\epsilon}
\def\e{{\rm e}}
\def\sgn{{\rm sgn}}

\theoremstyle{plain} 

\newtheorem{theorem}{Theorem}[section]
\newtheorem{lemma}{Lemma}[section]

\newtheorem{corollary}[theorem]{Corollary} \theoremstyle{definition}
\newtheorem{definition}{Definition}[section] \theoremstyle{remark}
\newtheorem{remark}{Remark}[section]

\newcommand{\R}{{\mathbb R}} 

\newcommand{\Hc}{{\mathcal H}}

\newcommand{\Z}{{\mathbb Z}}
\newcommand{\E}{{\mathcal E}}
\newcommand{\F}{{\mathcal F}}

\newcommand{\Tr}{{\mathcal T}}
\newcommand{\Ph}{{\mathcal P}}

\newcommand{\resto}{{\mathcal R}}

\newcommand{\T}{\mathbb{T}} 

\newcommand{\Poi}[2]{\left\{#1,#2\right\}}

\newcommand{\Or}{{\mathcal O}}

\numberwithin{equation}{section}

\begin{document}

\author{D. Bambusi, A. Carati, T. Penati\\ Universit\`a degli Studi di
Milano, \\ Dipartimento di Matematica ``F. Enriques'', \\ Via Saldini
50, 20133 Milano, Italy}

\title{Boundary effects on the dynamics of chains of coupled
oscillators}

\date{\today}

\maketitle

\begin{abstract}
We study the dynamics of a chain of coupled particles subjected to a
restoring force (Klein-Gordon lattice) in the cases of either periodic
or Dirichlet boundary conditions. Precisely, we prove that, when the
initial data are of small amplitude and have long wavelength, the main
part of the solution is interpolated by a solution of the nonlinear
Schr\"odinger equation, which in turn has the property that its
Fourier coefficients decay exponentially.  The first order correction
to the solution has Fourier coefficients that decay exponentially in
the periodic case, but only as a power in the Dirichlet case. In
particular our result allows one to explain the numerical computations
of the paper \cite{BMP07}.\end{abstract}

%\tableofcontentssults of some

%\setcounter{section}{0}

\section{Introduction.}
\label{S.0}

The dynamics of chains of coupled particles has been the object of a
huge number of studies, but only recently some numerical works (see
\cite{BG07, BMP07}) have shown that the boundary conditions have some
relevance on FPU type investigations. The goal of the present paper is
to study analytically the effects of the boundary conditions (BC) on
the dynamics of a simple 1-dimensional model, namely the so called
Klein Gordon lattice (coupled particles subjected to an on site
restoring force).  Precisely, we concentrate on the cases of periodic
and of Dirichlet boundary conditions, and use the methods of normal
form to study the dynamics. This leads to a quite clear understanding
of the role of the boundary conditions and to an explanation of the
numerical results of \cite{BMP07}.  On the contrary, our theory does
not allow one to explain the results of \cite{BG07}.

More precisely, we study the dynamics of a large lattice corresponding
to small amplitude initial data with long wavelength; we show that if
the size $N$ is large enough and the amplitude $\epsilon$ of the
initial excitation is of order $\mu:=\frac1{N}$, then the solution $z$
has the form
\begin{equation}
\label{rap.1}
z = \mu z^a(t)+\mu^2 z_{1}(t)\ ,
\end{equation}
up to times $|t|\leq\mathcal{O}(\mu^{-2})$. In \eqref{rap.1} $z^a$ is
interpolated by a solution of the nonlinear Schr\"odinger equation
(NLS) and has a behaviour which is independent of the BC. On the
contrary $z_{1}$ depends on the BC. Precisely, its Fourier
coefficients decrease exponentially in the periodic case, but only as
$|k|^{-3}$ in the Dirichlet case. 

The theory we develop in order to give the representation
\eqref{rap.1} provides a clear interpretation of the
phenomenon. Indeed, it turns out that the normal form of the system is
independent of the BC (and coincides with the NLS), but the coordinate
transformation introducing the normal form has properties which are
different in the periodic and in the Dirichlet case. In particular, in
the Dirichlet case it maps sequences which decay fast into sequences
which decay as $|k|^{-3}$. This introduces the slow decay in the
Dirichlet case.

It should be pointed out that our result still depends on the size $N$
of the lattice.\footnote{However our normal form holds in the region
of the phase points with small energy density $\epsilon^2$,
independently of N. The limitation $\epsilon\sim \mu$ comes from the
fact that we are only able to study the dynamics of the NLS in this
situation.} Nevertheless, we think that (within the range of validity
of our result) we clearly show the role of the boundary conditions and
provide a good interpretation of the numerical results.

\vskip 10pt 

The present situation has many similarities with the one occurring in
the theory of the Navier Stokes equation (see e.g. \cite{Tem90}),
where it is well known that the spectrum of the solution depends on
the boundary conditions. Moreover, we recall that a power law decay of
localized object has been previously observed in nonlinear lattice
dynamics in \cite{DP03,Pey04,Fla98,GF05} and that the connection
between the nonlinear Schr\"odinger equation and the dynamics of long
chains of particles was studied in many papers (see e.g. \cite{Kal89,
KSM92, Sch98}).

\vskip10pt

The paper is organized as follows: in sect. \ref{S.1} we present the
model, state our main result and discuss its relation with numerical
computations. In sect. \ref{S.3} we give the proof of the normal form
construction. In sect. \ref{S.4} we give the proof of the
decomposition \eqref{rap.1}. Some technical details are deferred to
the appendix. Each section is split into several subsections.

\vskip10pt

\noindent {\it Acknowledgments.} {This work was partially supported
  by MIUR under the project COFIN2005 ``Sistemi dinamici nonlineari e
  applicazioni fisiche''. We thank Antonio Ponno, Luigi Galgani and
  Simone Paleari for many discussion and suggestions during the
  preparation of this paper. We also thank Sergej Flach and Michel
  Peyrard for some interesting discussion.}

\newpage

\section{Main result}
\label{S.1}

In this chapter we present the model, we recall
some numerical simulations (see \cite{BMP07}) which clearly show the
dependence of the metastable Fourier decay on the boundary conditions
and we finally state the main theoretical result and use it to explain
the numerics.

\subsection{The model.}\label{model}

We consider a chain of particles described by the Hamiltonian function
\begin{eqnarray}
\label{mod}
H({p},{q}) &=& \sum_{j}\frac{p_j^2}{2}+\sum_j V(q_j)+\sum_j
 W(q_j-q_{j-1}) 
\end{eqnarray}
where $j$ runs from $0$ to $N$ in the case of Dirichlet boundary
conditions (DBC), namely $q_0=q_{N+1}=0$, while it runs from
$-(N+1)$ to $N$ in the case of periodic boundary conditions (PBC), i.e.
$q_{-N-1}=q_{N+1}$. The corresponding Hamilton equations are 
\begin{equation}
\label{eq.gen}
\ddot q_j=-V'(q_j)-W'(q_j-q_{j-1})+W'(q_{j+1}-q_{j})\ .
\end{equation} 
We recall that the standard Fermi Pasta Ulam model is obtained by
taking $V\equiv 0$ and $W(x)=\frac{x^2}{2}+\alpha
\frac{x^3}{3}+\beta\frac{x^4}{4}$. Here instead we will take
\begin{equation}
\label{mod.1}
V(x)=\frac{1}{2}x^2+\frac{1}{3}\alpha x^3+ \frac{1}{4}\beta x^4\
,\quad W(x)=\frac{1}{2}ax^2\ ,\ \alpha,\beta,a\geq0\ .
\end{equation}
Explicitly our Hamiltonian has the form
\begin{eqnarray}
\label{1.Ham} 
H &=& H_0+H_1+H_2, \\
\label{1.Ham1}
H_0(p,q) &:=& \sum_{j}\frac{p_j^2+q_j^2}{2}+a \frac{(q_j-q_{j-1})^2}2,
\\ 
\label{1.Ham2} \quad H_1(p,q) &:=& \alpha\sum_{j} \frac{q_j^3}3
 \ ,\quad H_2(p,q) := \beta \sum_{j}\frac{q_j^4}4\ .
\end{eqnarray}

\begin{remark}
\label{r.2}
In the case where
\begin{equation}
\label{sym}
V(x)=V(-x)
\end{equation}
the equations \eqref{eq.gen} with PBC are invariant under the
involution $q_j\mapsto -q_{-j}$, $p_j\mapsto -p_{-j}$. As a
consequence the submanifold of the periodic sequences which are also
skew-symmetric, is invariant under the dynamics. For this reason when
\eqref{sym} is fulfilled the case of DBC is just a subcase of the case
of PBC.  This happens in the standard FPU model and also in the case
of the Hamiltonian \eqref{mod} with the potential \eqref{mod.1} and
$\alpha=0$. The case $\alpha\not=0$ is the simplest one where a
difference between DBC and PBC is possible.
\end{remark}

Consider the vectors
\begin{equation}
\label{2.basis}
{\hat e}_{k}(j)=
\begin{cases}
\frac{\delta_{PD}}{\sqrt{N+1}}\sin{\left(\frac{jk\pi}{N+1}\right)},\qquad
k=1,\ldots,N,\cr
\frac1{\sqrt{N+1}}\cos{\left(\frac{jk\pi}{N+1}\right)},\qquad
k=-1,\ldots,-N,\cr
\frac1{\sqrt{2N+2}},\qquad\qquad\qquad k=0,\cr
\frac{(-1)^j}{\sqrt{2N+2}},\qquad\qquad\qquad k=-N-1,\cr
\end{cases}
\end{equation}
then the Fourier basis is formed by $\hat e_k$, $k=1,\ldots,N$ and
$\delta_{PD}=\sqrt2$ in the case of DBC, and by $\hat e_k$,
$k=-N-1,\ldots,N$ and $\delta_{PD}=1$ in the case of PBC. Here we will
treat in a unified way both the cases of DBC and PBC, thus we will not
specify the set where the indexes $j$ and $k$ vary. Introducing the
rescaled Fourier variables $(\hat p_k,\hat q_k)$ defined by
\begin{eqnarray}
\label{2.nm}
p_{ j}=\sum_{ k}{\sqrt{\omega_{ k}}{\hat p}_{ k}{\hat e}_{ k}({
  j})},\qquad q_{ j}=\sum_{ k}{\frac{{\hat q}_{ k}}{\sqrt{\omega_{
  k}}}{\hat e}_{ k}({ j})},
\end{eqnarray}
where the frequencies are defined by
\begin{equation}
\label{2.freq}
\omega_{k} = \sqrt{1+4a {\sin^2{\left(\frac{k\pi}{2N+2}\right)}}},
\end{equation}
the Hamiltonian $H_0$ is changed to
\begin{equation}
\label{2.startH}
H_0=\sum_{ k}{\omega_{ k}\frac{{\hat p}^2_{ k}+{\hat q}^2_{ k}}2}\ .
\end{equation}

\subsection{The phenomenon and its numerical evidence.}
\label{phenomenon}

\begin{figure}[t]
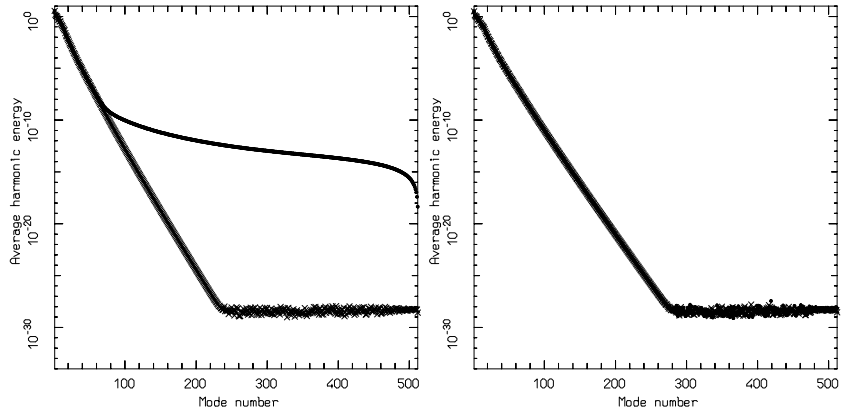

\begin{center}
\includegraphics[width=.45 \textwidth]{FIG1}
\includegraphics[width=.45 \textwidth]{FIG2}
\end{center}
\caption{Averaged harmonic energies distribution. DBC (dots) and PBC
  (crosses) with $N=511,\,a=.5,\,\E = 0.001,\,T=10^5$. Panel (a):
  $\alpha$ model with $\alpha=0.25$. Panel (b): $\beta$ model with
  $\beta=0.25$.}
\label{fig2}
\end{figure}

\begin{figure}[t]
\begin{center}
\includegraphics[angle=-90, width=.47 \textwidth]{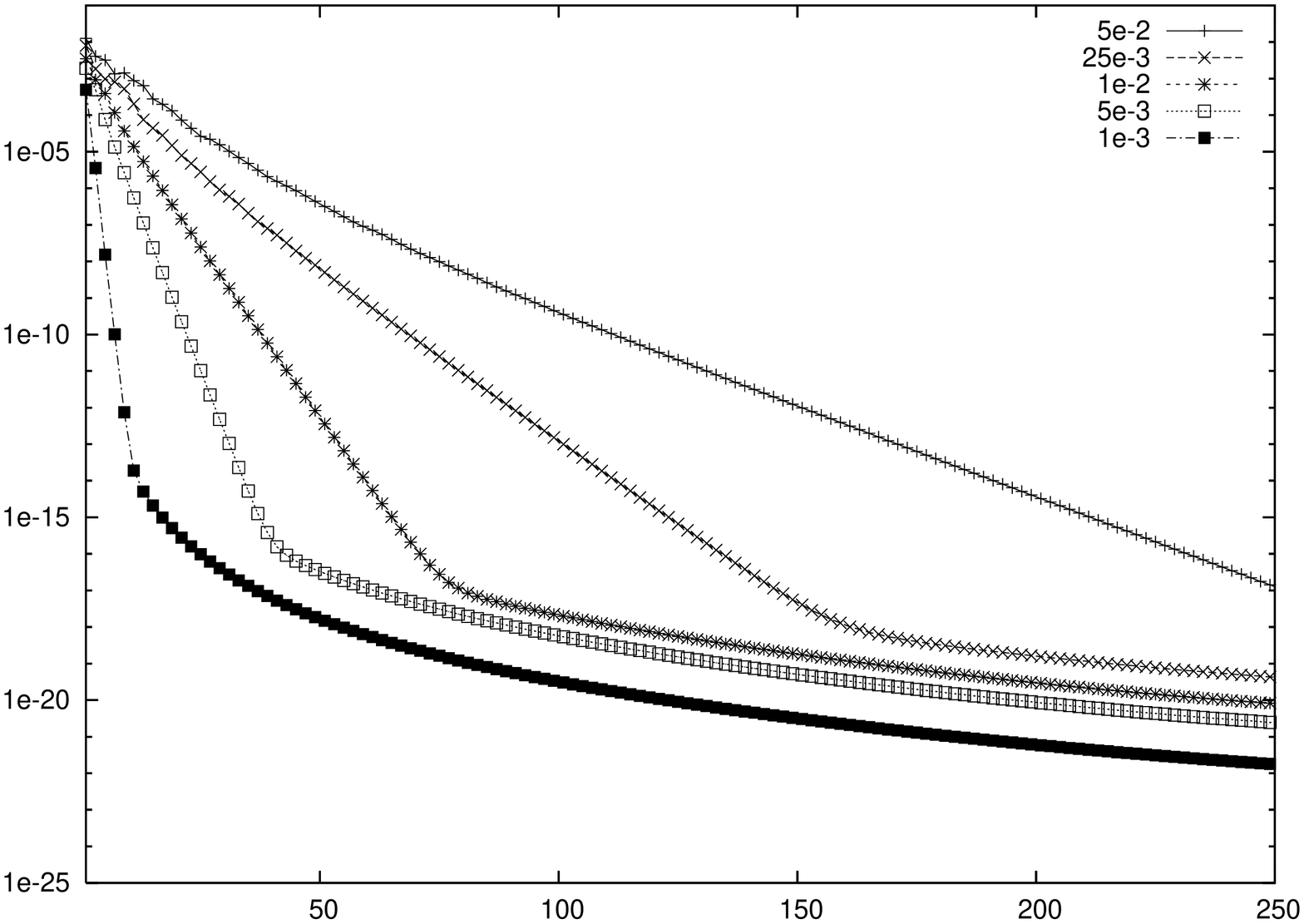}
\includegraphics[angle=-90, width=.47 \textwidth]{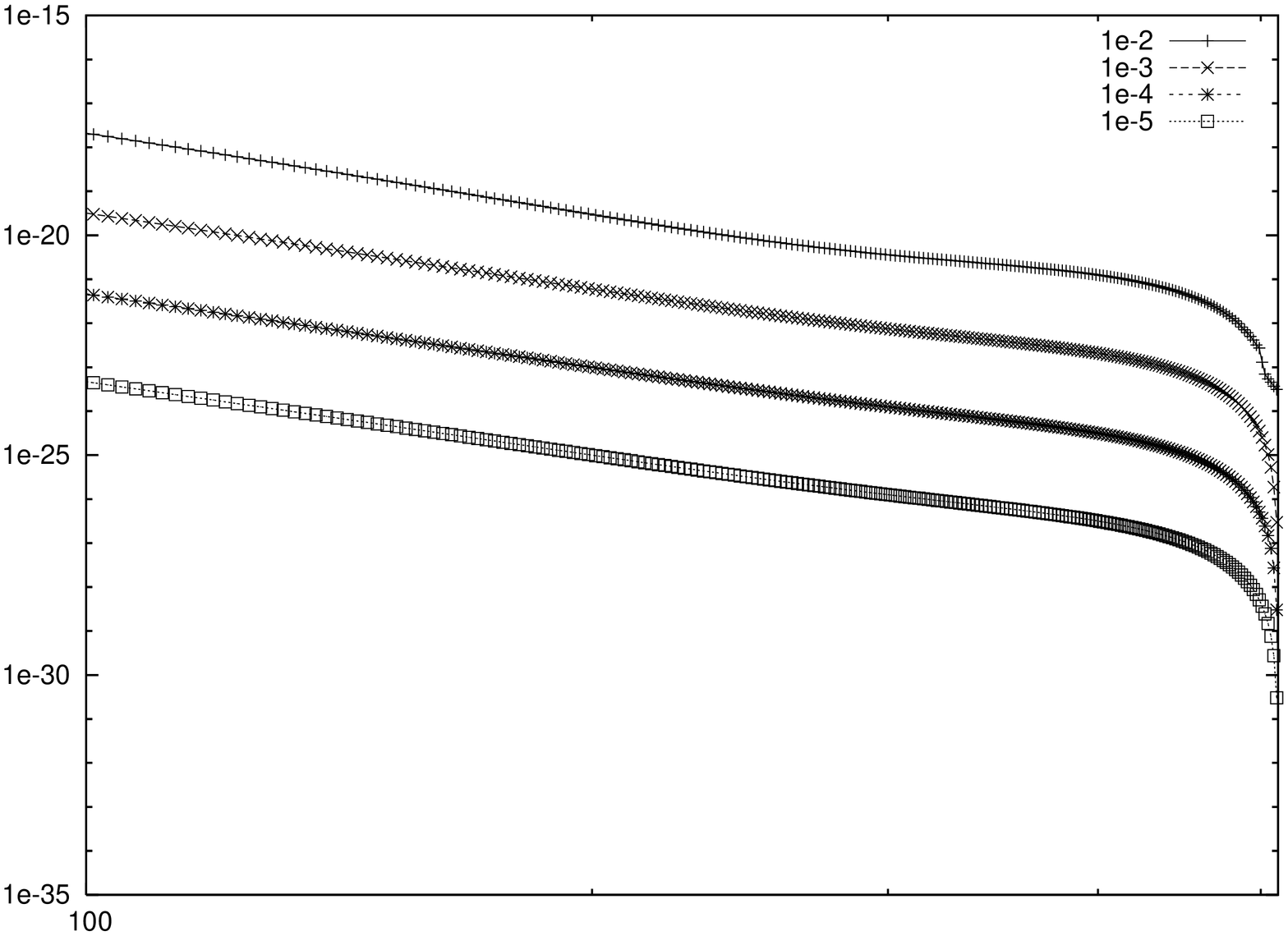}
\end{center}
\caption{DBC with parameters
    $N=511,\,a=0.5,\,\alpha=0.1,\,T=10^5$. Panel (a): distribution of
    ${\langle E_k\rangle}$ in semi-log scale. Energy densities:
    $\E=0.05,0.025,0.01,0.005,0.001$. Panel (b): distribution of
    ${\langle E_k\rangle}$ in log-log scale. Energy densities:
    $\E=0.01,0.001,0.0001,0.00001$.}
\label{fig3}
\end{figure}

Let us define the energy of a normal mode and its time average by
\begin{displaymath}
E_k:=\omega_k\frac{\hat p_k^2+\hat q_k^2}{2},\qquad\qquad \langle
E_k\rangle(t) := \frac{1}{t} \int_0^tE_k(s)ds\ .
\end{displaymath}
In the case of PBC the oscillators of index $k$ and $-k$ are in
resonance, so the relevant quantity to be observed is the average
$\overline{\langle E_k\rangle}= \frac12 (\langle E_k\rangle +\langle
E_{-k}\rangle) $.

Take an initial datum with all the energy concentrated on the first
Fourier mode with energy density $\E\equiv H_0/N = 0.001$.
Integrating the system numerically one can see that after a short
transient the time averages of the harmonic energies relax to well
defined steady values, which persist for very long times.  In figures
\ref{fig2} we plot in a semi-log scale the time-average energies
$\langle E_k\rangle (T)$ (or $\overline {\langle E_k\rangle}(t)$) at
time $T= 10^5$ (subsequent to the relaxation time) as a function of
the index $k$. The parameters in the two panels are
$\alpha=0.25,\,\beta=0$ and $\alpha=0,\,\beta=0.25$ respectively. In
both the cases $a=0.5$. In each distribution the dots refer to the DBC
case while the crosses pertain the PBC one\footnote{Actually we plot
only modes with odd index, since, as shown in \cite{BMP07} the
dynamics involves only modes with odd index t}.

While in panel (b) one clearly observes a perfect overlapping of the
exponential part of the decays, in panel (a) a sharp difference
arises. Indeed, while the PBC solution is once more characterized by
an exponential distribution, in the case of DBC one sees a richer
behavior: at an energy approximately equal to $10^{-8}$ there is
crossover and a new regime appears. Nevertheless, a striking
similarity among the exponential part of the two dynamics is evident.

To describe more carefully the situation in the case of DBC we plot in
figures \ref{fig3} four different distributions of the quantities
${\langle E_k\rangle}$, in a semi-log and in a log-log scale
respectively. They correspond to different values of the energy
density (see the caption). In the first panel we plot the first part
of the distribution: we notice that by decreasing the energy density,
the slope of the exponential decay of the low frequencies
increases. In the second figure, instead, we focus our attention on
the second part of the distribution: we see that the corresponding
curves are parallel. So a change of energy only induces a
translation. Except for the last part, that we will interpret as due
to discreteness effects, the curves are very well interpolated by a
straight line giving a power decay with an exponent close to $-6$. A
similar behavior is also obtained if one excites a few modes of large
wave length.

\subsection{Explanation of the phenomenon.}
\label{explanation}

In order to state our main result we need a topology in the phase
space.

\begin{definition}
\label{phasespace}
Let us define the spaces $\ell^2_{s,\sigma}$ of the sequences
$p=\{\hat p_k\}$ s.t.
\begin{equation}
\label{e.1.4}
\norma{p}^2_{s,\sigma} := \mu\sum_{k}[k]^{2s}{\rm
  e}^{2\sigma|k|}\left|\hat p_k \right|^2<\infty\ ,\qquad
      [k]:=\max\left\{1,|k|\right\}\ ,
\end{equation}
and the phase spaces $\Ph_{s,\sigma}:=\ell^2_{s,\sigma}
\times\ell^2_{s,\sigma}\ni( p,q)$. 
\end{definition}

The main part of the solution will be described by NLS, so we consider
a a smooth solution $\varphi(x,\tau)$ of the nonlinear Schr\"odinger
equation
\begin{equation}
\label{nls.111}
-\im\partial_\tau\varphi = -\partial_{xx}\varphi +
\gamma\varphi|\varphi|^2\ ,\quad
\gamma:=\frac{3}{8a}\left(\beta-\frac{10}{9}\alpha^2 \right)\ ,\quad
x\in\T:=\R/2\pi\Z
\end{equation}
For fixed $\tau$ we will measure the size of a function $\varphi$ by
the norm
\begin{equation}
\label{e.13.5}
\norma \varphi_{s,\sigma}^2:=\sum_k[k]^{2s}\e^{2\sigma|k|}|\hat
\varphi_k|^2
\end{equation}
where $\hat \varphi_k$ are the Fourier coefficients of $\varphi$, which are
defined by 
\begin{equation}
\label{inifi}
\varphi(x):=\sum_{k}\hat \varphi_k\hat
e^c_k(x)\ ;
\end{equation}
here 
$\hat e^c_k(x)$ is the continuous Fourier basis,
\begin{equation}
\label{he.300}
\hat e^c_k:=\left\{
\begin{matrix}
\frac{1}{\sqrt\pi}\cos kx & k>0 \\ \frac{1}{\sqrt{2\pi}} & k=0 \\
\frac{\delta_{PD}}{\sqrt\pi}\sin(-kx) & k<0
\end{matrix}
\right.\ .
\end{equation}
\begin{remark}
\label{r.e}
The dynamics of \eqref{nls.111} is well known \cite{FadTak,GKa03}. 
Precisely, if
$\tilde\gamma <\tilde\gamma^*$ then $\forall \sigma\geq 0$ there
exists $0\leq\sigma'<\sigma$ such that from
$\norma{\varphi(x,0)}_{s,\sigma}=1$ it follows
$\norma{\varphi^a(x,t)}_{s,\sigma'}\leq C_s\sim 1$ for all times.
\end{remark}

Corresponding to $\varphi(x,\tau)$ we define an approximate solution
of the original model by
\begin{equation}
\label{z.app}
z^a_j(t)\equiv\left(p_j^a(t),q_j^a(t)\right):=\left(\Re\left(e^{it}\varphi(\mu
j,a \mu^2 t)\right)\ ,\ \Im\left(e^{it}\varphi(\mu j,a \mu^2 t)
\right) \right)
\end{equation}

Our main result concerns the comparison between $z^a(t)$ and the
solution $z(t)$ of the original system with initial datum
\begin{displaymath}
z(0):= \mu z^a(0)
\end{displaymath}

\begin{theorem}
\label{t.22}
Assume $\norma{\varphi^0(x)}_{s,\sigma}\leq1$, then $\forall T>0$ there
exists $\mu^*>0$ with the following properties: if $\mu<\mu_*$ then
there exists $z_1(t)$ defined for $|t|\leq T\mu^{-2}$ such
that
\begin{equation}
\label{3.gf}
z(t)=\mu z^a(t)+\mu ^2z_1(t)
\end{equation}
where $z^a$ is the approximate solution just defined and,
\begin{equation}
\label{t.2.1}
\norma{z_1(t)}_{s,\sigma'}\leq C
\end{equation}
with 
\begin{equation}
\label{lsa1}
\begin{matrix}
\frac{1}{2}<s\ ,\quad&\sigma'> 0&\ \text{if\ PBC}
\\
\frac12< s<\frac{5}{2}\ ,&\sigma'=0&\  \text{if\ DBC}
\end{matrix}
\end{equation}
\end{theorem}
The above result gives an upper estimate of the error $z_1(t)$. We
want now to compute it, at least for short times. To this end we
assume, for simplicity, that $\varphi^0$ is purely imaginary,
i.e. that the initial datum has zero velocity. Consider now the
complex function $\psi_{10}$ with Fourier coefficients given by
\begin{equation}
\label{phi.1}
(\psi_{10})^{\wedge}_k:=\frac{\alpha}{6\sqrt2}\quadr{4e^{\im\omega_kt} -
  3e^{2\im t} - e^{-2\im t} - 6\im e^{\im\omega_kt} + 6\im} \hat\Phi_k
\end{equation}
where $\Phi(x):=-\im \varphi^0(x)^2$, and let
$z_{10}\equiv(p_{10},q_{10})$ 
be the sequence with 
$$ \frac{p_{10,j}+\im q_{10,j}}{\sqrt2}=\psi_{10}(\mu j)\ .
$$
We have the following
\begin{theorem}
\label{t.22.1}
For any $0<b<1$,
one has 
\begin{equation}
\label{lsb}
\norma{z_1(t)-z_{10}(t)}_{s,\sigma'}\leq C\mu^{1-b}\ ,\quad |t|\leq
\frac{T}{\mu^{b}}
\end{equation}
\end{theorem}

The above theorems provide the interpretation for the numerical
results of the previous paragraph. Referring to FIG.~1, we identify
the exponential part of the distribution as due to the main part of
the solution, namely $z^a$, which is the same for both the boundary
conditions. In the case of PBC, the inequality \eqref{t.2.1} implies
that also the error $z_1$ is exponentially decreasing; thus the whole
solution is in particular exponentially localized in Fourier space.

In the case of DBC the situation is different. Indeed the correction,
namely $z_1$, is ensured to have coefficients such that the series
with general term $|k|^{2s} |\hat z_{1,k}(t)|^2$ is convergent; this
is very close to say that
\begin{displaymath}
  |\hat z_{1,k}(t)|^2< \frac{C}{|k|^{2s+1}}
\end{displaymath}
which, taking $s$ very close to $5/2$ essentially gives a power law
decay like $|k|^{-6}$. 

Then theorem \ref{t.22.1} shows that this is actually optimal, as seen
by taking 
\begin{equation}
\label{ne.1}
\varphi^0(x)=\im \sin x
\end{equation}
(as in the numerical computations). Indeed in such a case one has
\begin{equation}
\label{ne.2}
\hat \Phi_k=\left(\sin^2 x\right)^{\wedge}_k\sim\frac{1}{k^3}
\end{equation}
so that \eqref{lsb} shows that after a time of order 1 the energy of
the $k$-th mode is of order $\mu^2/k^6$ as shown by the numerics. 

\jumpsec

%%%%%%%%%%%%%%%%%%%%%%%%%%%%%% NORMAL FORM CONSTRUCTION

\section{The normal form construction.}
\label{S.3}

In this part we introduce and use the methods of the normal form
 theory for the proof of our main result.  Accordingly one looks for a
 canonical transformation putting the system in a simpler form.

\subsection{Preliminaries and main claim.}
\label{prel}
We first need to introduce some notations:

\begin{itemize}
\item $z$ will denote a phase point. In particular a phase point can
be represented using the coordinates $(p_j,q_j)$ of the lattice's
particles or the Fourier coordinates $(\hat p_k, \hat q_k)$.

\item In the phase space we will also use coordinates $\psi$ defined
  by
$$
\hat\psi_k:=\frac{\hat p_k+\im\hat q_k}{\sqrt2}
$$
and, in real space
$$ 
\psi_j:=\sum_{ k}{\hat \psi}_{ k}{\hat e}_{ k}({ j})
$$

\item Given a Hamiltonian function $H$, we will denote by $X_H$ the
corresponding Hamiltonian vector field. Thus if one uses for example
the variables $(p_j,q_j)$, one has
$$ X_H(p,q)=\left(-\frac{\partial H}{\partial q_j},\frac{\partial
H}{\partial p_j}\right)\ \quad\text{or}\quad X_H(\psi):=\left(\im
\frac{\partial H}{\partial \bar \psi_j} \right)
$$ Correspondingly we will write the Hamilton equations of a
Hamiltonian function $H$ by
\begin{displaymath}
\dot z = X_H(z)\ .
\end{displaymath}

\item The Lie transform $\Phi^1_\chi$ generated by a Hamiltonian function
 $\chi$ is the time one flow of the corresponding Hamilton equations,
 namely
\begin{displaymath}
\Phi^1_\chi:=\Phi^t_\chi\big|_{t=1}\ ,\quad
 \frac{d}{dt}\Phi^t_{\chi}(z) =
 X_{\chi}\tond{\Phi^t_{\chi}(z)}\ ,\quad
 \Phi^t_\chi\big|_{t=0}=\uno\ .
\end{displaymath}

\item The Poisson bracket $\poisson{f}{g}$ of two function $f,g$
is defined by
\begin{displaymath}
 \poisson{f}{g}:=d fX_g\equiv \sum_{j}\left( \frac{\partial
f}{\partial q_j} \frac{\partial g}{\partial p_j} -\frac{\partial
f}{\partial p_j} \frac{\partial g}{\partial q_j}\right).
\end{displaymath}

\end{itemize}

The normalizing transformation will be constructed by composing two
Lie transforms $\Tr_1=\Phi^1_{\chi_1}$ and $\Tr_2=\Phi^1_{\chi_2}$
generated by two functions $\chi_1$ and $\chi_2$. Taking $\chi_1$ and
$\chi_2$ to be homogeneous polynomials of degree 3 and 4 respectively,
an elementary computations shows that
\begin{eqnarray}
\label{e.1.1a}
H\circ\Tr_1\circ\Tr_2&=& H_0
\\
\label{e.1.1b}
&+&\left\{\chi_1,H_0\right\}+H_1
\\
\label{e.1.1c}
&+&\left\{\chi_2,H_0\right\}+\frac{1}{2}\left\{\chi_1,\left\{\chi_1,
H_0\right\}\right\}+\left\{\chi_1,H_1\right\}+H_2
\\
&+&{\rm h.o.t.}
\end{eqnarray}
where the term \eqref{e.1.1b} is a homogeneous polynomial of degree 3,
the term \eqref{e.1.1c} has degree 4 and h.o.t denotes higher order
terms. We will construct a function $\chi_1$ such that \eqref{e.1.1b}
vanishes and we will show that there exists a $\chi_2$ such that
\eqref{e.1.1c} is reduced to normal form.

To make precise the construction we need to split $H_0$ as follows
\begin{eqnarray}
\label{splh0}
H_0 &=& H_{00}+H_{01}\ ,\qquad H_{00}(z) := \sum_{k}\frac{\hat
  p_k^2+\hat q_k^2}{2}, \\ 
\nonumber H_{01}(z) &:=&
  \sum_{k}\nu_k\frac{\hat p_k^2+\hat q_k^2}{2} \ ,\qquad
  \nu_k:=\omega_k-1=
  a\frac{4{\sin^2{\left(\frac{k\pi}{2N+2}\right)}}}{\omega_k+1}\leq 2a
\end{eqnarray}

\begin{definition}
\label{def.nor.for}
A polynomial $Z$ will be said to be in normal form if it Poisson
commutes with $H_{00}$, i.e. if
\begin{equation}
\label{nor.for.1}
\Poi{H_{00}}{Z}\equiv 0\ .
\end{equation}
\end{definition}

\begin{remark}
\label{sub}
In order to study the system with DBC we will {\it always} extend the
system to a system defined for $j=-(N+1),...,N$ with PBC, which is
invariant under the involution $q_j\mapsto -q_{-j}$, $p_j\mapsto
-p_{-j}$. The extension of the $\beta$-model is obtained without
modifying the equations, while the extension of the $\alpha$-$\beta$
model is given by the system
\begin{equation}
\label{e.34}
\dot q_j=p_j\ ,\quad \dot p_j=-q_j-a(\Delta_1q)_j - \alpha s_j q_j^2 -
\beta q_j^3,
\end{equation}
where $s_j$ is a discretization of the step function given by
\begin{equation}
\label{C.e.10}
s_j:=\left\{ 
\begin{matrix}1&\text{if}\ j\geq 1
\\
0&\text{if}\ j=0
\\
-1&\text{if}\ j\leq 1
   \end{matrix}
  \right.
\end{equation}
and 
\begin{equation}
\label{2.dlap}
(\Delta_1 q)_{j}=2q_{j}-q_{j+1}-q_{j-1},
\end{equation}
is the discrete Laplacian. {\em The need of the introduction of the
sequence $s_j$ is at the origin of the finite smoothness of the
solution in the DBC case}.
\end{remark}

\smjump

We are going to prove the following
\begin{theorem}
\label{main}
Assume $a<1/3$. Then, there exists an analytic canonical
transformation $z=\Tr(\psi)$, defined in a neighborhood of the origin
\begin{equation}
\label{6.17}
H\circ \Tr= H_0(\psi)+Z(\psi)+\resto(\psi)\ ,
\end{equation}
where $Z$ is in normal form and the following holds true
\begin{itemize}
\item[1)] The remainder, the normal form and the canonical
  transformation are estimated by
\begin{eqnarray}
\label{1.R)}
\norma{X_{\resto}(\psi)}_{s,\sigma} &\leq& C_\resto \norma{\psi}_{s,\sigma}^4,\\
\norma{X_{Z}(\psi)}_{s,\sigma} &\leq& C_Z\norma{\psi}_{s,\sigma}^3\ ,\\
\label{2)}
\norma{\psi-\Tr(\psi)}_{s,\sigma} &\leq& C_\Tr \norma{\psi}_{s,\sigma}^2\
\end{eqnarray}
where 
\begin{equation}
\label{lsa}
\begin{matrix}
\frac{1}{2}<s\ ,\quad&\sigma\geq 0&\ \text{if\ PBC}
\\
\frac12< s<\frac{5}{2}\ ,&\sigma=0&\  \text{if\ DBC}
\end{matrix}
\end{equation}
and $C_{\resto,Z,\Tr}$ are constants independent of $N$.
\item[2)] One has $\Tr=\uno+X_{\chi_{01}}+\resto_{\Tr}$ with
$\chi_{01}$ given by
\begin{equation}
\label{3)}
\chi_{01}(\psi)=-\frac\alpha{6\sqrt{2}}\sum_j{\left(\frac13\psi_j^3 -
  \frac\im3\overline\psi_j^3 + 3\psi_j\overline\psi_j^2 -
  3\im\psi_j^2\overline\psi_j\right)}.
\end{equation}  
and 
\begin{equation}
\label{3.chir}
\norma{\resto_{\Tr}(\psi)}_{s-s_1,\sigma}\leq
C_{\resto_{\Tr}}\mu^{s_1} \norma{\psi}_{s,\sigma}^2
\end{equation}
where the parameters vary in the range 
\begin{equation}
\label{range.1}
\begin{matrix}
0\leq s_1< s-\frac12<2 ,\ \sigma=0\ \text{for\ DBC} \\ 0\leq
s_1<s-\frac12\ ,s_1\leq 2\ ,\ \sigma\geq 0\ \text{for\ PBC}
\end{matrix}
\end{equation}

\item[3)] The normal form has the following structure 
\begin{equation}
\nonumber Z(\psi,\bar \psi)=Z_0+Z_r\
\end{equation}
where 
\begin{equation}
Z_0:= \tilde\gamma\sum_j\left|\psi_j\right|^4 ,\quad
\tilde\gamma:=\frac{3}{8}\left(\beta-\frac{10}{9}\alpha^2 \right)
\end{equation}
and
\begin{equation}
\label{4.r}
\norma{X_{Z_r}(\psi)}_{s-s_1,\sigma}\leq \mu^{s_1} C_r
\norma{\psi}_{s,\sigma}^3\ 
\end{equation}
with the parameters varying in the range \eqref{range.1}.
\end{itemize}
\end{theorem}

\smjump

The proof of this theorem is divided into three parts. In the first
one we will prove an abstract normal form theorem under the assumption
that the non-linearity corresponding to the system \eqref{1.Ham} is
smooth. In the second part we will prove such a smoothness
assumption. In the third part we will compute the first order term of
the normal form and of the transformation and we will estimate the
corresponding errors.

We point out that there are 3 delicate points in the proof: the first
one is to solve the homological equation (see lemma \ref{sol.homo});
the second one is to prove smoothness of the perturbation in the
optimal space, and the third one is the actual computation of the main
part of the normal form and of the canonical transformation.
\smjump

\subsection{An abstract normal for theorem.}
\label{abs}

First we recall that a homogeneous polynomial map
$F:\Ph_{s,\sigma}\to\Ph_{s,\sigma }$ of degree $r$ is continuous and
also analytic if and only if it is bounded, i.e. if there exists a
constant $C$ such that
\begin{equation}
\label{p.e.1a}
\norma{F(z)}_{s,\sigma}\leq C\norma{z}^r_{s,\sigma}\qquad\qquad \forall
z\in\Ph_{s,\sigma}.
\end{equation}

\begin{definition}
\label{poly}
The best constant such that \eqref{p.e.1a} holds is called the norm of
$F$, and will be denoted by $\opnorma F_{s,\sigma}$. One has
\begin{equation}
\label{poly.1}
\opnorma F_{s,\sigma}:=\sup_{\norma
z_{s,\sigma}=1}\norma{F(z)}_{s,\sigma}\ .
\end{equation}
\end{definition} 

\begin{definition}
\label{p.d.2}
A polynomial function $f$, homogeneous of degree $r+2$, will be said
to be of class $\poly_{s,\sigma}^r$ if its Hamiltonian vector field
$X_f$ is bounded as a map from $\Ph_{s,\sigma}$ to itself.
\end{definition}

\begin{theorem}
\label{abnf}
Let $H_0$ be as above (see \eqref{2.startH}) with $a<1/3$. Assume that
$H_j\in\poly_{s,\sigma}^{j}$ for some fixed $s,\sigma$ and for
$j=1,2$; then there exists an analytic canonical transformation
$\Tr=\Tr_1\circ\Tr_2$, defined in a neighborhood of the origin in
$\Ph_{s,\sigma}$ such that
\begin{equation}
H\circ \Tr= H_0+Z+\resto\ ,
\end{equation}
where $Z$ is in normal form and the remainder, the normal form and the
canonical transformation are bounded by
\begin{eqnarray}
\label{1.ab.R}
\norma{X_{\resto}(z)}_{s,\sigma} &\leq& C_\resto \norma{z}_{s,\sigma}^4\ ,\\ 
\norma{X_{Z}(z)}_{s,\sigma} &\leq& C_Z\norma{z}_{s,\sigma}^3\ ,\\
\norma{z-\Tr(z)}_{s,\sigma} &\leq& C_\Tr \norma{z}_{s,\sigma}^2,\ 
\end{eqnarray}
with constants $C_{\resto,Z,\Tr}$ depending only on
$a,\opnorma{X_{H_1}}_{s,\sigma},\opnorma{X_{H_2}}_{s,\sigma}$.
\end{theorem}

The rest of this subsection will be occupied by the proof of theorem
\ref{abnf}. First we need some simple estimates. 

\begin{lemma}
\label{p.c.1}
Let $f\in\poly_{s,\sigma}^r$ and $g\in\poly_{s,\sigma}^{r_1}$, then
$\left\{f,g\right\}\in\poly_{s,\sigma}^{r+r_1}$, and 
\begin{equation}
\label{p.c.2}
\opnorma{X_{\{f;g \}}}_{s,\sigma}\leq(r+r_1+2)
\opnorma{X_f}_{s,\sigma} \opnorma{X_g}_{s,\sigma}
\end{equation}
\end{lemma}
\proof First remember that
\begin{equation}
\label{p.c.3}
X_{\Poi{f}{g}} = [X_f;X_g]=dX_f X_g - dX_g X_f\ .
\end{equation}
We recall now that, given a polynomial $X$ of degree $r+1$, there
exists a unique $(r+1)-$linear symmetric form $\widetilde X$ such that
\begin{displaymath}
X(z)=\widetilde X(z,...,z)\ ,
\end{displaymath}
then, \eqref{p.c.3} is explicitly given by
\begin{equation}
\label{p.c.44}
X_{\{f;g\}}=(r+1)\widetilde X_f(X_g(z),z,z,...,z)-(r_1+1)\widetilde
X_g(X_f(z),z,z,...,z)
\end{equation}
moreover from \eqref{poly.1} one has
\begin{equation}
\label{p.c.4}
\norma{\widetilde X(z_1,...,z_{r+1}) }_{s,\sigma}\leq
\opnorma{X}_{s,\sigma}
\norma{z_1}_{s,\sigma}... \norma{z_{r+1}}_{s,\sigma}
\end{equation}
from which the thesis immediately follows. \qed

\begin{remark}
\label{p.r.2a}
Let $f\in\poly_{s,\sigma}^r$ then the corresponding vector field
generates a flow in a neighborhood of the origin in
$\Ph_{s,\sigma}$. 
\end{remark}

\begin{lemma}
\label{p.l.2}
Let $\chi$ be of class $\poly_{s,\sigma}^{r}$, and let
$f\in\poly^{r_1}_{s,\sigma}$. Let $\Phi^1_\chi$ be the Lie transform
generated by $\chi$, then each term of the Taylor expansion of
$f\circ\Phi^1_\chi$ is a polynomial with bounded vector field.
\end{lemma}

\proof Iterating the relation 
\begin{displaymath}
\frac{d}{dt}(f\circ\Phi^t_\chi) = \Poi{\chi}{f}\circ\Phi^t_\chi 
\end{displaymath}
one gets that the Taylor expansion of $f\circ\Phi_\chi^1$ is given
by
\begin{eqnarray*}
f\circ \Phi^1_\chi &=& \sum_{l\geq 0}{f_l}\\
f_0&=&f,\qquad\qquad f_l=\frac1l\Poi{\chi}{f_{l-1}},\qquad l\geq1.
\end{eqnarray*}
Then, the thesis follows from lemma \ref{p.c.1}\qed

A key role in the proof of theorem \ref{abnf} is played by the so
called homological equation, namely
\begin{equation}
\label{homo}
\left\{\chi_j;H_0\right\}+f_j=Z_j
\end{equation}
where $f_j\in\poly_{s,\sigma}^j$ is a given polynomial, and $\chi_j\in
\poly_{s,\sigma}^j$, $Z_j\in\poly_{s,\sigma}^j$ are to be determined
with the property that $Z_j$ is in normal form.

\begin{lemma}
\label{sol.homo}
Consider the Homological equation \eqref{homo} with $f_j$ of class
$\poly_{s,\sigma}^{j}$, with $j=1,2$. Assume that $a<\frac{1}{3}$, then
\eqref{homo} admits a solution $\chi_j$, $Z_j\in
\poly_{s,\sigma}^{j}$ with
\begin{equation}
\label{sti.homo}
\opnorma{X_{\chi_j}}_{s,\sigma} \leq
\frac1{2(1-3a)}\opnorma{X_{f_j}}_{s,\sigma}.
\end{equation}
\end{lemma}

\proof First we rewrite the homological equation as 
\begin{equation}
\label{hom.1}
(L_0+L_1)\chi_j=f_j-Z_j
\end{equation}
where the operators $L_0$ and $L_1$ are defined by
$$
L_0\chi_j:=\Poi{H_{00}}{\chi_j}\ ,\quad L_1\chi_j:=\Poi{H_{01}}{\chi_j}
$$ and $H_{00}$ and $H_{01}$ are defined by \eqref{splh0}. We will
invert $L_0$ and solve \eqref{hom.1} by Neumann series (see
\cite{BDGS07}).

We begin by showing that the space $\poly_{s,\sigma}^j$, $j\leq2$
decompose into the sum of the kernel ${\rm Ker}(L_0)$ of $L_0$ and of
its range ${\rm Im}(L_0)$. Moreover, we show that $L_0$ is invertible
on its range.

Given $f\in \poly_{s,\sigma}^j $ with $j=1,2$ define
\begin{eqnarray}
\label{2.int}
Z&:=& \frac1T\int_0^T{f\left(\Psi^t(z)\right)dt}\ ,
\\
\label{hom.4}
\chi&:=&\frac1T\int_0^T{\quadr{t
f\left(\Psi^t(z)\right)-Z\left(\Psi^t(z)\right)}dt}\ ,
\end{eqnarray}
where $\Psi^t$ is the flow of $X_{H_{00}}$ and $T=1$ is its
period. Then an explicit computation shows that $Z\in$Ker$(L_0)$, and
that (see \cite{BG93})
\begin{equation}
\label{homo.aux}
L_0\chi=f-Z\ .
\end{equation}
Thus denoting by $Q$ the projector on the kernel of $L_0$, and
$P=\uno-Q$ the projector on the range, one sees that the \eqref{2.int}
is a concrete definition of $Q$, while \eqref{hom.4} is the definition
of $L_0^{-1}$ restricted to Im($L_0$). It remains to show that
$Z,\chi\in\poly_{s,\sigma}^j$. Remark that, since $\Psi^t$ is a
canonical transformation one has
\begin{eqnarray}
\label{hom.5}
X_Z(z)&\equiv& X_{Qf} (z)= \frac1T\int_0^T \tond{\Psi^{-t}\circ X_
{f}\circ\Psi^t}(z)d t\ , \\ \label{hom.6} X_\chi(z)&\equiv&
X_{L_0^{-1}Pf}(z) = \frac1T\int_0^T \tond{\Psi^{-t}\circ X_
{Pf}\circ\Psi^t}(z)tdt\ ,
\end{eqnarray}
From which it follows that
\begin{eqnarray}
\label{hom.7}
\opnorma{X_{Qf}}_{s,\sigma}&\leq& \opnorma{X_f}_{s,\sigma}\ ,
\\
\label{hom.8}
\opnorma{X_{ {L_0^{-1}Pf }}}_{s,\sigma}&\leq&\frac{1}{2}
\opnorma{X_{Pf}}_{s,\sigma}\leq\frac{1}{2}
\opnorma{X_{f}}_{s,\sigma}\ \Longrightarrow \norma{L_0^{-1}}\leq
\frac{1}{2}
\end{eqnarray}
where the last norm is the norm of $L_0^{-1}$ as a linear operator
acting on the space $\poly_{s,\sigma}^j$, and thus
$Qf\in\poly_{s,\sigma}^j$, $L_0^{-1}Pf\in\poly_{s,\sigma}^j$. 

We come now to the true homological equation \eqref{hom.1}. We look
for a solution $\chi_j=P\chi_j$ and $Z_j=QZ_j$. Applying $P$ or $Q$ to
\eqref{hom.1}, remarking that since $[L_0,L_1]=0$ one has
$[P,L_1]=[Q,L_1]=0$, we get 
\begin{equation}
\label{hom.2}
(L_0+L_1)\chi_j=Pf_j\ ,\quad Qf_j=Z_j\ .
\end{equation}
The first of \eqref{hom.2} is formally solved by Neumann
series, i.e. defining
\begin{equation}
\label{hom.3}
(L_0+L_1)^{-1}:=\sum_{k\geq0}(-1)^k\left(L_0^{-1}L_1\right)^kL_0^{-1}\
,\quad  \chi_j:=(L_0+L_1)^{-1}Pf_j\ . 
\end{equation}
To show the convergence of the series in operator norm we need an
estimate of $\norma{L_1}$. To this end remark that, for any $s,\sigma$
one has
\begin{displaymath}
\opnorma{X_{H_{01}}}_{s\sigma}\leq 2a\ , 
\end{displaymath}
which using lemma \ref{p.c.1} implies $\norma{L_1}\leq 2a(j+1)\leq
6a$. It follows that the series \eqref{hom.3} converges provided
$a<1/3$, which is our assumption, and that
\begin{equation}
\label{sti.L0L1}
\norma{(L_0+L_1)^{-1}} \leq \frac1{2(1-3a)},
\end{equation}
which concludes the proof.
\qed

\noindent 
{\it End of the proof of Theorem \ref{abnf}}. From lemma
\ref{sol.homo} one has that the solution $\chi_1$ of the homological
equation with $f_1\equiv H_1$ is well defined {\it provided
$H_1\in\poly_{s,\sigma}^1$ for some $s,\sigma$}. Then $\chi_1$
generates a Lie transform $\Tr_1$ which puts the system in normal form
up to order 4. Then the part of degree four of $H\circ\Tr_1$ takes the
form
\begin{equation}
\label{e1.}
f_2:=\frac{1}{2}\left\{\chi_1,\left\{\chi_1,
H_0\right\}\right\}+\left\{\chi_1,H_1\right\}+H_2\equiv\frac12
\left\{\chi_1,H_1\right\}+H_2
\end{equation}
which is of class $\poly_{s,\sigma}^2$. It follows that one can use
the homological equation with such a known term and determine a
$\chi_2$ which generates the Lie transformation putting the system in
normal form up to order 4. This concludes the proof of theorem
\ref{abnf}.\qed

\subsection{Proof of the smoothness properties of the nonlinearity}
\label{smnon}

In this subsection we prove the following lemma

\begin{lemma}
\label{p.l.3}
Let $H_j$, $j=1,2$ be given by \eqref{1.Ham}. Consider the vector
fields $X_{H_j}$ of the cubic and of the quartic terms of the
Hamiltonian: they fulfill the estimates
\begin{equation}
\label{p.f.1}
\norma{\mod{X_{H_1}}(z)}_{s,\sigma}\leq G_1 \norma{z}_{s,\sigma}^{2} \
,\quad \left\{
\begin{matrix}
\frac{1}{2}<s, &\sigma\geq 0& \text{if\ PBC}
\\
\frac12<s<\frac{5}{2}, &\sigma=0& \text{if\ DBC}
\end{matrix}
\right.
\end{equation}
and 
\begin{equation}
\label{p.f.2}
\norma{\mod{X_{H_2}}(z)}_{s,\sigma}\leq G_2 \norma{z}_{s,\sigma}^{3}\
,\quad \frac{1}{2}<s\ ,\ \sigma\geq 0\,\ \text{both\ cases}
\end{equation}
where we set $G_j:=\opnorma{X_{H_j}}_{s,\sigma}$.
\end{lemma}

The proof will be split into two parts. First we show that it is
possible to prove the result working on interpolating functions, and
then we show that the ``interpolating nonlinearities'' have a smooth
vector field when the parameters $s,\sigma$ vary in the considered
range.
\begin{remark}
\label{re.non}
Define $T_l$ as the map $q_j\mapsto q_j^{l+1}$ in the case of PBC, and
$[T_1(q)]_j=q_j^2s_j$ and $[T_2(q)]_j:= q_j^{3}$ in the case of
DBC. Then the vector field $X_{H_l}$ has only $p$ components, moreover
the norms are defined in terms of the Fourier variables, so we have to
estimate the map constructed as follows
\begin{eqnarray*}
&&\hat q_k\mathop{\mapsto}^1 \hat q_k^S := \frac{\hat
q_k}{\sqrt{\omega_k}}\mathop{\mapsto}^{\F} q_j:=\sum_{k} \hat q_k^S
{\hat e}_{ k}({ j})\mathop{\mapsto} T_l(q_j)\mathop{\mapsto}
^{\F^{-1}} \\ 
&&
\mathop{\mapsto} ^{\F^{-1}} p_k^S := \sum_j q_j^{l+1} {\hat
e}_{ k}({ j})\mathop{\mapsto}^2 \hat p_k:=\sqrt{\omega_k} p_k^S
\end{eqnarray*}
It is immediate to realize that the maps 1 and 2
are smooth (the frequencies are between 1 and 3) so it is enough to
estimate the remaining maps. In turns the remaining maps essentially
coincide with the map $T_l(q)$ read in terms of standard Fourier
variables (without the factors $\sqrt{\omega_k}$). These are the maps we
will estimate.
\end{remark}

{\it All along this section we will use a definition of the Fourier
coefficients of a sequence not including the factors
$\sqrt{\omega_k}$, namely we define $\hat q_k$ by
$$
q_j=\sum_{k}\hat q_k\hat e_k(j)\ .
$$
}

We start by showing how to use the interpolation in order to make
estimates. To this end we define an interpolation operator $I$ by
\begin{equation}
\label{e.13.3}
[I(q)](x):=\sum_k \sqrt\mu\hat q_k \hat e_k^c(x)
\end{equation}

We also define a restriction operator $R$ that to a function
associates the corresponding sequence, by
\begin{equation}
\label{e.13.4}
[R(u)]_j:=u(\mu j)
\end{equation}
We remark that the operator $R$ is defined on functions which do not
necessarily have finitely many non-vanishing Fourier coefficients. 

\begin{remark}
\label{r.33.3}
With the definition \eqref{e.13.5} one has 
\begin{equation}
\label{e.13.7}
\norma{Iq}_{s,\sigma}=\norma{q}_{s,\sigma}\ .
\end{equation}
\end{remark}
\begin{lemma}
\label{smo.1}
For any $s>1/2$ there exists a constant $C_6(s)$ such that one has
\begin{equation}
\label{smo.2}
\norma{Ru}_{s,\sigma}\leq C_6\norma{u}_{s,\sigma}
\end{equation}
\end{lemma}
\proof Denote $q_j=(Ru)_j$; using the formula
\begin{equation}
\label{e.11.2}
\hat e^c_k(\mu j)=\hat e^c_{k+2(N+1)m}(\mu j)=\frac1{\sqrt{\mu}}\hat e_k(j)
\end{equation}
one gets
\begin{displaymath}
q_j=\sum_{k\in\Z}\hat u_k\hat e^c_k(\mu
j)=\sum_{k=-(N+1)}^N{\sqrt\mu}{\hat e^c_k(\mu j)}
\sum_{m\in\Z}\frac{\hat u_{k+2(N+1)m}}{\sqrt\mu} =
\sum_{k=-(N+1)}^N{\hat e_k(j)\hat q_k}\ ,
\end{displaymath}
from which 
\begin{equation}
\label{smo.4}
\hat q_k=\sum_{m\in\Z}\frac{\hat u_{k+2(N+1)m}}{\sqrt\mu}.
\end{equation}
Let's define
\begin{displaymath}
\gamma_{m,k}(s,\sigma):=\frac{[k]^{s} e^{\sigma |k|}}{[k+2(N+1)m]^{s}
  e^{\sigma |k+2(N+1)m|}} =
\tond{\frac{[k]}{[k+2(N+1)m]}}^s\frac{e^{\sigma |k|}}{e^{\sigma
    |k+2(N+1)m|}}
\end{displaymath}
and replace \eqref{smo.4} in the norm $\norma{q}_{s,\sigma}^2$, then
we get
\begin{eqnarray*}
\norma{q}_{s,\sigma}^2 &=&\mu\sum_{k=-(N+1)}^N{[k]^{2s} e^{2\sigma
    |k|}|\hat q_k|^2}=\\ &=&\sum_{k=-(N+1)}^N[k]^{2s} e^{2\sigma
  |k|}\left|\sum_{m\in\Z} \hat u_{k+2(N+1)m} \right|^2
=\\ &=&\sum_{k=-(N+1)}^N[k]^{2s} e^{2\sigma |k|}\left|\sum_{m\in\Z}
\frac{\gamma_{m,k}(s,\sigma)}{\gamma_{m,k}(s,\sigma)}\hat
u_{k+2(N+1)m} \right|^2\leq\\ &\leq&
\sum_{k=-(N+1)}^N\left(\sum_{m\in\Z}{\gamma_{m,k}^2(s,\sigma)}\right)\sum_{m\in\Z}
    {\frac{[k]^{2s} e^{2\sigma |k|}}{\gamma_{m,k}^2(s,\sigma)}}|\hat
    u_{k+2(N+1)m}|^2\leq\\ &\leq& C_6(s)\norma{u}^2_{s,\sigma}\ .
\end{eqnarray*}
Indeed, since $k=-(N+1),\ldots,N$, we can estimates the two factors of
$\gamma_{m,k}$ as follows:
\begin{itemize}

\item 
\begin{displaymath}
|2(N+1)m+k|\geq|2(N+1)|m|-|k||\geq
\begin{cases}
|k|,\qquad m=0,\\
N+1\geq|k|,\qquad m\not=0,
\end{cases}
\end{displaymath}
which gives
\begin{displaymath}
\frac{e^{\sigma |k|}}{e^{\sigma |k+2(N+1)m|}}\leq 1.
\end{displaymath}
\item
\begin{displaymath}
\frac{[k]}{[k+2(N+1)m]} = 
\begin{cases}
\frac1{[2(N+1)m]}<\frac1{[2m]},\qquad\qquad m=0,\\
\frac{|k|}{|k+2(N+1)m|}\leq \frac1{1+2m}<\frac1{[2m]},\qquad m\not=0
\end{cases}
\end{displaymath}
which gives
\begin{displaymath}
\tond{\frac{[k]}{[k+2(N+1)m]}}^{2s}\leq\frac1{[2m]^{2s}}.
\end{displaymath}
\end{itemize}
\qed

\begin{corollary}
\label{l.13.1}
Let $T:\R^{2(N+1)}\to \R^{2(N+1)}$ be a polynomial map, assume that
there exists an ``interpolating polynomial map $T^c$'' such that
$T=RT^cI$. If the map $T^c$ is bounded in some space
$H^{s,\sigma}$, with $s>1/2$, 
then $T$ is bounded in
$\ell^2_{s,\sigma}$. Moreover one has
\begin{equation}
\label{e.13.15}
\opnorma{\mod T}_{s,\sigma}\leq C_6\opnorma{\mod {T^c}}_{s,\sigma}\ .
\end{equation}
\end{corollary}

First we define the interpolating maps we have to study. They are
$T^c_l(u):=u^{l+1}$ in the case of PBC and $T^c_1(u)=\sgn(x) u^2(x)$ and 
$T^c_2(u):=u^{3}$ in the case of DBC. Here we introduced the function 
$$\sgn(x):=\left\{ 
\begin{matrix}1&\text{if}\ x>0
\\
0&\text{if}\ x=0
\\
-1&\text{if}\ x<1
   \end{matrix}
  \right.
$$ The estimates \eqref{p.f.1} {in the case of PBC} and \eqref{p.f.2}
for both boundary conditions are proved in Lemma \ref{A.2} by a
standard argument on the Sobolev norm of the product of two
functions. We come to the estimate of $T^c_1$ in the case of DBC.

We will denote by $H^s_o $ the subspace of $H^{s,0}$ composed by the
odd functions $u(x)$ on $[-\pi,\pi]$.

\begin{lemma}
\label{A.3}
For any $1/2<s<5/2$, The operator
\begin{displaymath}
T^c_1(u):=u^2 \mathrm{sgn}(x)
\end{displaymath}
is smooth from $H^s_o$ in itself and there exists $C_7(s)$ such that
\begin{displaymath}
\norma{T^c_1u}_{s}\leq C_7\norma{u}^2_{s}.
\end{displaymath}
\end{lemma}

\proof We begin with the case $2\leq s <\frac52$. First, observe that
the function $T^c_1(u)$ is odd when $u$ is odd. we will prove the
thesis by showing that the second weak derivative $d^2T^c_1(u)$ of
$T^c_1(u)$ is in $H_o^{s-2}$. First remark that, by an explicit
computation which exploits the fact that $u(0)=0$ one has
\begin{eqnarray*}
d^2T_1^c(u) = 2[\sgn(x)(u'^2)+\sgn(x)ud^2u]\ .
\end{eqnarray*}
We show now that both terms are in $H^{s-2}_o$. The second term can be
considered as the product of the function $d^2u\in H^{s-2}_o$ and of
$\sgn(x)u$. This last function is of class $H^1$, as it is seen by
computing its derivative, namely
$$
d\left(\sgn(x)u\right)=  \sgn(x)du(x) +
\delta(x)u(x)=\sgn(x)du(x)
$$ which clearly belongs to $L^2$. From lemma \ref{A.2} it follows
that the product $\sgn(x)u[d^2u]\in H^{s-2}$.

Concerning the term $\sgn(x)(u'^2)$, it can be considered as the
product of $u'^2\in H^{s-1}$ and of $\sgn(x)$, which is of class
$H^r_o$ for all $r<1/2$, as it can be seen by explicitly computing
its Fourier coefficients. Thus lemma \ref{A.2} gives the result.  

The case $1\leq s<2$ is easier and works in a very similar
way. Indeed, since $u(x)\sgn(x)\in H^1$ and $du(x)\in H^s$ with $0\leq
s<1$, the derivative $d(u^2(x)\sgn(x)) = 2u(x)(du(x))\sgn(x)$ belongs
to $H^s$ with $0\leq s<1$, which gives the thesis.

Quite different is the case $\frac12<s<1$, since by hypothesis no weak
derivative exists for the function $u(x)$. We exploit the
following equivalent definition of the
norm of the Sobolev space $H^s_o([-\pi,\pi])$ with real exponent $s$
\begin{equation}
\label{HSdef}
\norma{u}^2_s:=
\int_{-\pi}^\pi\int_{-\pi}^\pi{\frac{|u(x)-u(y)|^2}{|x-y|^{1+2s}}dx dy},
\end{equation}
and the skew-symmetry of the periodic function $u(x)\in H^s_o$. We
want to prove that $g(x):=u(x)\sgn(x)\in H^s$ with $\frac12<s<1$; more
precisely
\begin{equation}
\label{holder}
\norma{u(x)\sgn(x)}^2_s < 4\norma{u(x)}^2_s,\qquad\qquad \frac12<s<1.
\end{equation}
The symmetries of $g(x)$ on the given domain allow to simplify the
integral in \eqref{HSdef}
\begin{eqnarray*}
\norma{g}^2_s &=&
2\int_{0}^\pi\int_{-\pi}^\pi{\frac{|g(x)-g(y)|^2}{|x-y|^{1+2s}}dx
dy}=\\ &=&
2\int_{0}^\pi\int_{0}^\pi\quadr{{\frac{|g(x)-g(y)|^2}{|x-y|^{1+2s}}} +
\frac{|g(x)-g(y)|^2}{|x+y|^{1+2s}}}dxdy=\\ &=&
2\int_{0}^\pi\int_{0}^\pi\quadr{{\frac{|u(x)-u(y)|^2}{|x-y|^{1+2s}}} +
\frac{|u(x)-u(y)|^2}{|x+y|^{1+2s}}}dxdy\leq\\ &\leq&
4\int_0^\pi\int_0^\pi{\frac{|u(x)-u(y)|^2}{|x-y|^{1+2s}}}dxdy <
4\norma{u}^2_s.
\end{eqnarray*}

\qed

\subsection{Computation of the normal form and of the transformation.}
\label{s.4}

In this section we will concentrate on the case of DBC which is the
most difficult one.

Consider again the Hamiltonian \eqref{1.Ham}, introduce complex
variables $\xi_j,\eta_j$ defined by
\begin{equation}
\label{C.e.8}
\xi_j=\frac{p_j+\im q_j}{\sqrt2}\ ,\qquad \eta_j=\frac{p_j-\im
  q_j}{\im\sqrt2}
\end{equation}
and split $H_0 = \Hc_{00} + \Hc_{01}$ with
\begin{eqnarray}
\label{h00n}
\Hc_{00}&:=&\sum_{j}\frac{p_j^2+q_j^2}{2}\equiv \im \left\langle
\xi;\eta \right\rangle_{\ell^2}\ ,
\\\label{h01n}
\Hc_{01}&:=&\sum_{j}\frac{a}{2}(q_j-q_{j-1})^2\equiv
\frac{a}{2}\left\langle \frac{\xi-\im\eta}{\im\sqrt
2};-\Delta_1 \frac{\xi-\im\eta}{\im\sqrt 2} \right\rangle_{\ell^2} \ ,
\end{eqnarray}
where $\left\langle\xi;\eta\right\rangle_{\ell^2}:=\sum_j\xi_j\eta_j$.
\begin{remark}
\label{spli.1}
The above splitting is different from the one introduced in
\eqref{splh0} which had been used in the proof of lemma
\ref{sol.homo}, and which was based on the Fourier variables. In
particular one has $\left\{H_{00};H_{01}\right\}\equiv 0$, but
$\left\{\Hc_{00};\Hc_{01}\right\}\not\equiv 0$.
\end{remark}
In the variables $(\xi,\eta)$ the flow $\Phi^t$ of $\Hc_{00}$ acts as
follows
\begin{equation}
\label{flowH00}
\Phi^t(\xi,\eta) := 
\begin{cases}
\xi_j\mapsto e^{\im t}\xi_j,\\
\eta_j\mapsto e^{-\im t}\eta_j.
\end{cases}
\end{equation}
The third order part of the Hamiltonian takes the form
\begin{equation}
\label{C.e.9}
\Ham{0}{1}=\frac\alpha3\sum_j\left(\frac{\xi_j-\im\eta_j}{\im\sqrt 2}
\right)^3 s_j
\end{equation}
where $s_j$ is the discrete step function defined in
\eqref{C.e.10}. The form of $H_2$ will be given below.
 
Denoting
\begin{equation}
\label{C.e.2}
\lie_0:=\left\{ \Hc_{00},.\right\}\ ,\quad \lie_1:=\left\{ \Hc_{01};.\right\}
\end{equation}
we rewrite the homological equation for $\chi_1$ as follows
\begin{equation}
\label{ham.7}
(\lie_0+\lie_1)\chi_1=\Ham{0}{1}
\end{equation}
which is solvable since the Kernel of $\lie_0+\lie_1=L_0+L_1$ on
polynomials of third order is empty. The solution $\chi_1$ of
\eqref{ham.7} is unique and, as shown in Lemma \ref{sol.homo}, exists.

By a direct computation one has
\begin{equation}
\label{chi1.split}
\chi_1 = \lie_0^{-1}H_1 - (\lie_0+\lie_1)^{-1}\lie_1\lie_0^{-1}H_1;
\end{equation}
we are going to show that the second term is much smaller than the
first one. Before starting, a couple of remarks are in order.

\begin{remark}
\label{C.r.1a}
The discrete Laplacian is $\ell^2$-symmetric on periodic sequences
\begin{equation}
\label{l2sym}
<\Delta_1\xi,\eta>_{\ell^2}=<\xi,\Delta_1\eta>_{\ell^2}.
\end{equation}
This is an immediate consequence of the fact that in Fourier variables
it acts as a multiplier by a real factor.
\end{remark}

\begin{remark}
\label{C.r.1}
In Fourier coordinates the discrete Laplacian $\Delta_1$ defined in
\eqref{2.dlap} acts as a multiplier by $\sin^2 k\mu$. It follows that
it has norm 1 when acting on anyone of the spaces
$\ell^2_{s,\sigma}$. Moreover, since
\begin{displaymath}
|\sin^{2}(k\mu)|\leq k^{s_1}\mu^{s_1},\qquad s_1\in[0,2],\qquad
 k\mu\in\quadr{0,\pi}
\end{displaymath}
one also has
\begin{equation}
\label{C.e.3a}
\norma{\Delta_1 \xi}_{s-s_1,\sigma}\leq \mu^{s_1}\norma
\xi_{s,\sigma},\qquad s_1\in[0,2].
\end{equation}
From \eqref{C.e.3a} and \eqref{h01n} it follows
\begin{equation}
\label{H01}
\norma{X_{\Hc_{01}}(\xi,\eta)}_{s-s_1,0}\leq
C(a)\mu^{s_1}\norma{(\xi,\eta)}_{s,0},\qquad 0\leq s_1< s-\frac12.
\end{equation}
\end{remark}

\begin{lemma}
\label{C.l.1}
Assume $a<\frac13$, then $\chi_1=\chi_{10}+\chi_{1r}$ with
\begin{equation*}
\chi_{10}(\xi,\eta)=\frac\alpha{6\sqrt{2}}\sum_j{\left(\frac13\xi_j^3
  - \frac\im3\eta_j^3 + 3\xi_j\eta_j^2 - 3\im\xi_j^2\eta_j\right)},
\end{equation*}  
and there exists $C_8(a,G_1)$ such that
\begin{equation}
\label{c.l.2}
\norma{X_{\chi_{1r}}(\xi,\eta)}_{s-s_1,\sigma}\leq C_8\mu^{s_1}
\norma{(\xi,\eta)}_{s,\sigma}^2 ,
\end{equation}
with
\begin{displaymath}
\left\{
\begin{matrix}
0\leq s_1< s-\frac12<2 ,\ \sigma=0\ \text{for\ DBC}
\\ 0\leq s_1<s-\frac12,\,s_1\leq 2\, \ \sigma\geq 0\ \text{for\ PBC}
\end{matrix}
\right.
\end{displaymath}
\end{lemma} 

\proof According to \eqref{chi1.split}, let's define
\begin{displaymath}
\chi_{10}:=\lie_0^{-1} H_1,\qquad\qquad
\chi_{1r}:=-(\lie_0+\lie_1)^{-1}\lie_1\chi_{10}=-(L_0+L_1)^{-1}\lie_1\chi_{10}.
\end{displaymath}
Since $\chi_{10}$ solves the homological equation
$\Poi{\Hc_{00}}{\chi_{10}} = H_1$, it can be explicitly computed by
\begin{eqnarray*}
\chi_{10}(\xi,\eta) &=& \int_0^1{H_1\tond{\Phi^t(\xi,\eta)}dt}=\\ &=&
-\frac{\alpha}{6\sqrt{2}}\sum_j{s_j\left(-\frac13\xi_j^3 +
\frac\im3\eta_j^3 - 3\xi_j\eta_j^2 + 3\im\xi_j^2\eta_j\right)}
\end{eqnarray*}
which also implies
\begin{displaymath}
\opnorma{X_{\chi_{10}}}_{s,\sigma}\leq\frac12
\opnorma{X_{H_1}}_{s,\sigma}.
\end{displaymath}
The thesis follows from \eqref{p.c.44}, \eqref{sti.L0L1} and
\eqref{C.e.3a}.\qed

We move to the second Homological equation
\begin{equation}
\label{homo.2}
\Poi{\chi_2}{H_0} + \tilde H_2 = Z
\end{equation}
where $\tilde H_2 = H_2 + \frac12\Poi{\chi_1}{H_1}$ can be split
according to \eqref{chi1.split} into $\tilde H_2 = \tilde H_{20} +
\tilde H_{21}$ with
\begin{equation}
\label{f2.split}
\tilde H_{20}:= H_2 + \frac12\Poi{\chi_{10}}{H_1},\qquad \tilde H_{21}:=
\frac12\Poi{\chi_{1r}}{H_1}.
\end{equation}
More explicitly, the leading term $\tilde H_{20}$ is composed of
\begin{eqnarray*}
\frac12\Poi{\chi_{10}}{H_1} &=&
    \frac{\alpha^2}{24}\sum_j{\tond{\xi_j^4 + 4\im\xi_j^3\eta_j +
    10\xi_j^2\eta_j^2 - 4\im\xi_j\eta_j^3 + \eta_j^4}},\\ 
H_2 &=& \frac\beta{16}\sum_j{\tond{\xi_j^4 - 4\im\xi_j^3\eta_j -
    6\xi_j^2\eta_j^2 + 4\im\xi_j\eta_j^3 + \eta_j^4}}.
\end{eqnarray*}
Before proceeding, it is useful to perform the change of
variables 
\begin{equation}
\label{3.psi}
\hat\psi_k = \frac{\hat p_k+\im\hat q_k}{\sqrt2},\qquad\qquad \psi_j =
\sum_k{\hat\psi_k\hat e_k(j)}\ .
\end{equation}
which puts the quadratic part into diagonal form. This implies a
modification of $\tilde H_{20}$, which however is of higher order and
therefore will be included into the remainder terms. Indeed, the Lemma
below shows that the difference between $\xi$ and $\psi$ is small

\begin{lemma}
\label{3.l.1}
For any $2\geq s_1\geq 0$ it holds true
\begin{eqnarray}
\label{C.e.30}
\norma{\xi-\psi}_{s,\sigma} \leq
a\mu^{s_1}\norma{\psi}_{s+s_1,\sigma}\ ,\quad
\norma{\im\eta-\bar\psi}_{s,\sigma} \leq
a\mu^{s_1}\norma{\psi}_{s+s_1,\sigma}\ ,
\\
\nonumber
s>\frac12,\,\sigma\geq0.
\end{eqnarray}
\end{lemma}

\proof {By definition
\begin{eqnarray*}
\norma{\xi-\psi}_{s,\sigma}^2 &=& \sum_k{[k]^{2s}e^{2\sigma
    k}(\sqrt\omega_k - 1)^2(\frac{|\hat{\tilde q}_k|^2}{\omega_k} +
    |\hat{\tilde p}_k|^2)}\leq\\ &\leq&\sum_k{[k]^{2s}e^{2\sigma
    k}(\sqrt\omega_k - 1)^2(|\hat{\tilde q}_k|^2 + |\hat{\tilde
    p}_k|^2)};
\end{eqnarray*}
a Taylor expansion of the frequencies $\omega_k$ gives
\begin{displaymath}
|\sqrt\omega_k - 1|^2\leq \frac14
 a^2{4\sin^4{\left(\frac{k\pi}{2N+2}\right)}}\leq C
 a^2\mu^{2s_1}k^{2s_1}
\end{displaymath}
which is the thesis.

\qed
}
\begin{corollary}
\label{C.l.11}
In terms of the variables $\psi,\bar \psi$ one has $\tilde
H_2=H_{20}+H_{21}$ where
\begin{eqnarray}
\nonumber
H_{20}(\psi,\overline\psi) &:=&
    \frac{\alpha^2}{24}\sum_j{\tond{\psi_j^4 +
    4\im\psi_j^3\overline\psi_j - 10|\psi_j|^4 +
    4\im\psi_j\overline\psi_j^3 + \overline\psi_j^4}}+
\\
\label{c.l.r.1}
    &+&\frac\beta{16}\sum_j{\tond{\psi_j^4 -
    4\im\psi_j^3\overline\psi_j + 6|\psi_j|^4 -
    4\im\psi_j\overline\psi_j^3 + \overline\psi_j^4}},
\end{eqnarray}
and there exists $C_9(a,G_1)$ such that 
\begin{equation}
\label{c.l.r}
\norma{X_{H_{21}}(\psi,\bar \psi)}_{s-s_1,\sigma}\leq
C_9\mu^{s_1}\norma{\psi}_{s,\sigma}^3 
\end{equation}
\end{corollary}
Just averaging \eqref{c.l.r.1} with respect to the flow $\Phi^t$ it is
now immediate to get the following Corollary.

\begin{corollary}
\label{C.l.2}
The normal form $Z$ is composed of two terms, $Z=Z_0 + Z_r$, where the
leading term $Z_0$ is smooth and reads
\begin{equation}
\label{Z.split.1}
Z_0(\psi) =
\tilde\gamma\sum_j{|\psi_j|^4},\qquad\qquad\tilde\gamma:=\frac38\tond{\beta -
\frac{10}9\alpha^2}.
\end{equation}
while the remainder is small
\begin{equation}
\label{Z.split.2}
\norma{X_{Z_r}(\psi)}_{s-s_1,\sigma}\leq C_9\mu^{s_1}
\norma{\psi}_{s,\sigma}^3 \quad \left\{
\begin{matrix}
0\leq s_1<s-\frac12<2\ ,\ \sigma=0\ \text{for\ DBC} \\ 0\leq
s_1<s-\frac12,\,s_1\leq2,\, \sigma\geq 0\ \text{for\ PBC}
\end{matrix}
\right.
\end{equation}
\end{corollary}

Thus we have proved that the formula for $Z_0$ holds. The formula for
$\chi_{10}$ implies that the canonical transformation has the
structure \eqref{3)} and this concludes the proof.

\jumpsec

%%%%%%%%%%%%%%%%%%%%%%%%%%%%% NLS APPROXIMATION

\section{Proof of Theorem \ref{t.22}}
\label{S.4}

To discuss this issue we first write the equations of motion
of the first part of the normal form, namely of $H_0+Z_0$, in the form
\begin{equation}
\label{1.13}
\im\dot \psi_j=(A\psi)_{j}-\tilde\gamma\psi_j|\psi_j|^2\ ,
\end{equation}
where $A$ is a linear operator which in the Fourier variables acts as
a multiplier by $\omega_k=1+\frac{a}{2}\mu^2k^2+O(\mu^4k^4)$, namely
\begin{equation}
\label{1.14}
(\widehat{A\psi})_k=\omega_k\hat \psi_k =
(1+\frac{a}{2}\mu^2k^2)\psi_k +\mathcal{O}(\mu^4)\equiv
(\widehat{A_{NLS}\psi})_k + \mathcal{O}(\mu^4)\ .
\end{equation}
Take now an interpolating function for $\psi$, in other words a
function $u$ such that
\begin{equation}
\label{1.15}
\psi_j=\eps u(\mu j)\ ,
\end{equation}
where $\epsilon$ is a small parameter representing the
amplitude. Then, up to corrections of higher order, $u$ fulfills the
equation
\begin{equation}
\label{NLS}
-\im u_t=u + a\mu^2 u_{xx} - \tilde\gamma\eps^2 u|u|^2\ ,
\end{equation}
which, up to a Gauge transformation and a scaling of the time
introduced by
\begin{equation}
\label{twoscales}
u(x,t)= e^{\im t}\varphi(x,\tau)\ ,\qquad\qquad \tau := a\mu^2 t,
\end{equation}
gives the NLS equation
\begin{equation}
\label{NLS.1}
\im \varphi_\tau = - \varphi_{xx} + \gamma\varphi|\varphi|^2,
\qquad\qquad\gamma:=\frac{\tilde \gamma}{a}\frac{\eps^2}{\mu^2}.
\end{equation}

{\it In order to get a bounded value of $\gamma$, from now on we take
$\epsilon=\mu$}.

We now compare an approximate solution constructed through NLS and the
true solution of our Hamiltonian system. More explicitly,
corresponding to a solution $\varphi^a(t)$ of the NLS with analytic
initial datum, we define an approximate solution $\psi^a$ of the
original model by
\begin{equation}
\label{app.e}
\psi_j^a(t)=\mu e^{\im t}\varphi^a(\mu j,a\mu^2 t)\ .
\end{equation}
We also consider the true solution $\psi(t)$ of the Hamilton equation
of the original model, with initial datum $\psi_j^0 = \mu\varphi^a(\mu
j,0)$.

{\it From now on we will restrict to the case of DBC which is the
complicate one.}

We first work in the variables in which $H$ is reduced to the normal
form
\begin{equation}
\label{NF}
H = H_0 + Z_0 + Z_r + \resto.
\end{equation}

\begin{lemma}
\label{tilpsi1}
Let $\psi$ be the solution of the equations of motion of \eqref{NF}
with initial datum $\psi(0)=\Tr^{-1}(\psi^0)$ and let $\psi^a$ be as
defined in \eqref{app.e}, then $$\psi = \psi^a + \psi_1$$ with
\begin{displaymath}
\norma{\psi_1(t)}_{s,0}\leq C\mu^2,\quad \frac{3}{2}<s<\frac{5}{2}\
,\qquad |t|\leq \frac{T}{\mu^2}.
\end{displaymath}
\end{lemma}

\proof

Observe that the NLS equation for $\psi^a$ may be rewritten as
\begin{equation}
\label{psia.bis}
\dot \psi^a
=X_{H_{NLS,0}}(\psi^a)+X_{Z}(\psi^a)+X_\resto(\psi^a)-X_{\resto_1}(\psi^a)
\end{equation}
where
\begin{eqnarray*}
{\resto_1} &:=& \tond{{H_0}-{H_{NLS,0}}} + {Z_r} +
\resto\\ 
{H_{NLS,0}}(\hat\psi) &:=&
\sum_{k}\left(1+a\frac{1}{2}\mu^2k^2\right) \left|\hat
\psi_k\right|^2\ ,
\end{eqnarray*}
so that $X_{\resto_1}$ fulfills the estimate
\begin{eqnarray*}
\norma{X_{\resto_1}(\psi^a)}_{s,0} &\leq& C_1\mu^4
\norma{\psi^a}_{s+4,0} + C_r\mu^{s_1}
\norma{\psi^a}^3_{s+s_1,0} +
C_\resto\norma{\psi^a}^4_{s,0}\leq
\\ &\leq &C_4\mu^4\ ,\qquad\qquad \frac12< s<\frac{5}{2}\ .
\end{eqnarray*}
We compare $\psi^a$ with the full solution $\psi$ of the equation
\begin{equation*}
\dot\psi = X_{H_0}(\psi)+X_{Z}(\psi)+X_\resto(\psi)
\end{equation*}
with initial datum $\psi_0 = \Tr^{-1}(\psi^0)$, whose difference from
$\psi^0$ (initial datum for $\psi^a$) is controlled by
\begin{displaymath}
\norma{\psi_1(0)}_{s,0}=\norma{\psi^0-\Tr^{-1}(\psi^0)}_{s,0}\leq
C_\Tr\norma{\psi^0}_{s,0}^2\leq C_\Tr\mu^2.
\end{displaymath}
So we apply the Gronwall lemma (see lemma \ref{B.approx}) with
\begin{displaymath}
A:=X_{H_0}\ ,\qquad P=X_{Z}+X_{\resto}\ ,\qquad R= X_{\resto_1},
\end{displaymath}
obtaining that the error $\psi_1:=\psi-\psi^a$ from the NLS dynamics 
satisfies
\begin{equation}
\label{psi1.eq}
\dot\psi_1 = A\psi_1 + \quadr{P(\psi^a+\psi_1)-P(\psi^a)} +
X_{\resto_1}(\psi^a)
\end{equation}
and is estimated by
\begin{equation}
\label{psi1.est}
\norma{\psi_1(t)}_{s,\sigma'}\leq C_\Tr\mu^2 e^{C_6\mu^2
  t}+\frac{C_4\mu^4}{C_6\mu^2} \left(e^{C_6\mu^2 t}-1 \right)\leq
C_7\mu^2,\qquad |t|\leq T/\mu^2
\end{equation}
where $C_6:= 6C_Z$.

\qed

When we go back to the original variables, the solution $z$ may be
split as
\begin{displaymath}
z=\Tr(\psi) = \Tr(\psi^a+\psi_1) = \mu z^a + \mu^2 z_1 + \mu^3 z_2
\end{displaymath}
where we have defined
\begin{equation}
\label{psi1psi2}
z^a:=\frac{\psi^a}{\mu},\qquad z_1 := \frac{\psi_1 +
  X_{\chi_{01}}(\psi^a)}{\mu^2} ,\qquad \qquad
z_2 := \frac{z - \mu z^a - \mu^2z_1}{\mu^3}.
\end{equation}
More precisely, we claim that it holds
\begin{lemma}
\label{psi1}
We have
\begin{displaymath}
\norma{z_1}_{s,0}\leq C_1,\qquad\qquad
\norma{z_2}_{s,0}\leq C_2,
\end{displaymath}
up to times s.t. $|t|\leq\frac{T}{\mu^2}$. 
\end{lemma}

\proof The first inequality comes directly from the lemma
\ref{tilpsi1}. Concerning the second one, we remark that
\begin{eqnarray*}
z &=& \Tr(\psi^a+\psi_1) = \psi^a+\psi_1+\left[\Tr(\psi^a+\psi_1)-
  (\psi^a+\psi_1) \right]= \\ &=& \psi^a+\psi_1+
X_{\chi_1}(\psi^a+\psi_1) + \Or(\norma{\psi^a+\psi_1 }_{s,0}^3 )=
\\ &=& \psi^a+\psi_1+ X_{\chi_1}(\psi^a)+ \Or(\norma{\psi_1
}_{s,0}\norma{\psi^a }_{s,0}+ \norma{\psi^a+\psi_1 }_{s,0}^3 )
\end{eqnarray*}
(by differentiability of $X_{\chi_1}$ and Lagrange mean value
theorem). Finally from \eqref{c.l.2} with $s_1=1$ we have
$$
X_{\chi_1}(\psi^a)=X_{\chi_{01}}(\psi^a)+\Or(\mu\norma{\psi^a }^2_{s+1,0})
$$
\qed

We now analyze the first correction $z_1$. To this end we analyze
separately its two terms.

First remark that, from \eqref{psi1.eq} one has
$\psi_1=\psi_{10}+\Or(|t| \norma{\psi^a }^3_{s,0})$, where $\psi_{10}$
solves the equation $\dot \psi_{10}=A\psi_{10}$ with initial datum
$\Tr^{-1}(\psi^0)-\psi^0$. Thus we have
\begin{equation}
\label{e.333}
\psi_{10}=\e^{At}(
\Tr^{-1}(\psi^0)-\psi^0)=-\e^{At}X_{\chi_{01}}(\psi^0) +\Or(\mu
\norma{\psi^0 }^2_{s+1,0} )
\end{equation}
We now analyze the other term. To this end, with the aim of
considering the short time dynamics, we rewrite the equation
\eqref{psia.bis} as 
$$
\dot \psi^a=i\psi^a+\Or(\mu^3)\quad \iff\quad \psi^a=\e^{\im
  t}\psi^0+\Or(\mu^3|t|) 
\ .
$$
Thus exploiting the differentiability of $X_{\chi_{01}}$ we have
\begin{equation}
\label{e.334}
\mu^2z_1= -\e^{At}X_{\chi_{01}}(\psi^0)+X_{\chi_{01}}\left(\e^{\im
  t}\psi^0 \right)+ \Or(\mu^3|t|)\ .
\end{equation}
Since
\begin{eqnarray*}
\e^{At}X_{\chi_{01}}(\psi^0) &=&
-\frac\alpha{6\sqrt2}\quadr{3e^{At}\psi_0^2 + 6\im e^{At}|\psi_0|^2 +
  e^{At}\overline\psi_0^2},\\ X_{\chi_{01}}\left(\e^{\im t}\psi^0
\right) &=& -\frac\alpha{6\sqrt2} \quadr{3\psi_0^2e^{2\im t} + 6\im
  |\psi_0|^2 + e^{-2\im t}\overline\psi_0^2}\\
\end{eqnarray*}
which yields to
\begin{displaymath}
\mu^2 z_1 = -\frac\alpha{6\sqrt2}\quadr{3(e^{At}-e^{2\im t})\psi_0^2 +
  6\im (e^{At}-1)|\psi_0|^2 + (e^{At}-e^{-2\im t})\overline\psi_0^2} +
\Or(\mu^3|t|).
\end{displaymath}
In the case of $\psi_0 = \im\mu z_0$ (zero velocity initial datum) we
have
\begin{displaymath}
z_1 = \frac\alpha{6\sqrt2}\quadr{4e^{At} - 3e^{2\im t} - e^{-2\im t} -
6\im e^{At} + 6\im}z_0^2 + \Or(\mu|t|)
\end{displaymath}
which gives immediately the thesis.

\jumpsec
%%%%%%%%%%%%%%%%%%%%%%%%%%%%%%%%%%% APPENDIX A

\appendix
\section{Appendix: a few technical lemmas.}
\label{S.A.}

\begin{lemma}
\label{A.2}
Let $u\in H^{r,\sigma}$ and $v\in H^{s,\sigma}$ with $s>\frac12$ and
$s\geq r\geq0$, $\sigma\geq 0$. Then there exists $C=C(r,s)$ such that
the following inequality holds
\begin{equation}
\label{A.45}
\norma{uv}_{r,\sigma}\leq
C\norma{u}_{r,\sigma}\norma{v}_{s,\sigma}
\end{equation}
\end{lemma}

\proof In this proof it is useful to use the expansion of $u$ and $v$
on the complex exponentially. Thus we will write 
$u(x)=\sum_{k\in\Z}\hat u_k\e^{\im kx}/\sqrt{2\pi}$, and remark that
if in the definition of the norm cf \eqref{e.13.5} we substitute such
coefficients to the coefficients on the real Fourier basis, nothing
changes. This is due to the fact that both the basis of the complex
exponentials and the real Fourier basis are orthonormal. The advantage
is that in terms of the complex exponentials the product is mapped
into the convolution of the Fourier coefficients, thus we have simply
to estimate the norm of the function whose Fourier coefficients are 
\begin{equation}
\label{A.44}
(\hat u*\hat v)_k=\sum_{j} \hat u_{j-k}\hat v_k\ .
\end{equation}

As a preliminary fact we define the quantities 
\begin{displaymath}
\gamma_{j,k}=\frac{[j-k][k]^{\frac{s}{r}}}{[j]}\ ,
\end{displaymath}
and prove that there exists a constant $C(s,r)$ such that 
\begin{equation}
\label{A.12}
\sum_k{\frac1{\gamma_{j,k}^{2r}}}<C
\end{equation}
To obtain \eqref{A.12} we need some preliminary inequalities. For any
positive $a$ and $b$ one has 
\begin{equation}
\label{A.13}
(a+b)^{2r}\leq 2^{2r}\max\left\{a;b\right\}^{2r}< 4^r(a^{2r}+b^{2r})
\end{equation}
and for any $j$ and $k$ in $\Z$
\begin{equation}
\label{A.14}
\frac{1}{[j-k]^{2r}[k]^{2s-2r}}<
\frac1{(\min\{[j-k],[k]\})^{2s}}<\frac{1}{[j-k]^{2s}}+\frac{1}{[k]^{2s}}.
\end{equation}
From \eqref{A.13}, \eqref{A.14} and $[j] < [j-k]+[k]$ it follows
\begin{eqnarray*}
\frac1{\gamma_{j,k}^{2r}} &\leq&
\left(\frac{[j-k]+[k]}{[j-k][k]^{\frac{s}{r}}}\right)^{2r}\leq
\left(\frac1{[k]^{\frac{s}{r}}} +
\frac1{[j-k][k]^{\frac{s}{r}-1}}\right)^{2r}\leq\\ &\leq&
4^r\left(\frac1{[k]^{2s}} + \frac{1}{[j-k]^{2r}[k]^{2s-2r}}\right)\leq
\\ &\leq& 4^r \left(\frac1{[k]^{2s}} +
\frac{1}{[j-k]^{2s}}+\frac{1}{[k]^{2s}}\right)
\end{eqnarray*}
which gives \eqref{A.12} with
$$
C = 3\times 4^r \times \sum_{k\in\Z}\frac1{[k]^{2s}}.
$$
Hence 
\begin{eqnarray*}
\norma{uv}_{r,\sigma}^2 &=&\sum_{j}{[j]^{2r} e^{2\sigma
|j|}|\sum_k{\hat u_{j-k}\hat v_k}|^2}\leq\\ &\leq&
\sum_j{[j]^{2r}e^{2\sigma |j|}
\left(\sum_k{\frac1{\gamma_{j,k}^{2r}}}\right)
\left(\sum_k{\gamma_{j,k}^{2r}|\hat u_{j-k}\hat v_k|^2}\right)}\leq\\
&\leq& C^2\sum_{j,k}{[j-k]^{2r}e^{2\sigma |j-k|}|\hat
u_{k-j}|^2[k]^{2s}e^{2\sigma |k|}|\hat v_k|^2}\leq\\
&\leq& C^2 \tond{\sum_{l}{[l]^{2r}e^{2\sigma |l|}|\hat
u_{l}|^2}}\tond{\sum_{k}{[k]^{2r}e^{2\sigma |k|}|\hat
u_{k}|^2}},
\end{eqnarray*}
which concludes the proof.  \qed 

\vskip 24pt 

We state here a version of the Gronwall Lemma which is suited for our
estimates.  First we recall the following lemma.

\begin{lemma}
\label{B.Gronw}
Let $x:[0,T]\rightarrow \Ph$ be a differentiable function and $\Ph$ a
Banach space.  Assume that $\forall t\in[0,T]$ it fulfills the integral
inequality
\begin{equation}
\label{Gronw.hyp.1}
\norma{x(t)}\leq K+\int_0^t{\tond{a\norma{x(s)}+b}ds}
\end{equation}
with $a,b$ real and non negative parameters, then 
\begin{equation}
\label{gro.r}
\norma{x(t)}\leq e^{at}K+\frac{b}{a}\left(e^{at}-1\right)\ .
\end{equation}
\end{lemma}
The lemma we use in sect. \ref{S.4} is the following one.

\begin{lemma}
\label{B.approx}
Let $z(t),z^a(t)\in\Ph$, $t\in[-T,T]$ be respectively the solutions of
\begin{displaymath}
\begin{cases}
\dot z = Az + P(z),\\
z(0) = z_0
\end{cases}
,\qquad\qquad
\begin{cases}
\dot z^a = A z^a + P(z^a)
- R(z^a),\\
z^a(0) = z^a_0
\end{cases}
\end{displaymath}
where $A$ is the generator of a unitary group in $\Ph$, and
$\norma{z^a(t)}\leq C$. Assume also that the non-linearity $P$ has a
zero of third order at the origin and that for all
$t\in[-T,T]$ and all $z$ with $\norma{z}\leq 2C$
\begin{displaymath}
\norma{P(z)}\leq \rho_1\norma{z}^3,\qquad\qquad \norma{dP(z)}\leq
3\rho_1\norma{z}^2,
\end{displaymath}
and the remainder $R$ is estimated by
\begin{displaymath}
\norma{R(z^a(t))}\leq \rho_2,\qquad \forall t\in[-T,T]\ .
\end{displaymath}
Let $\delta:=z-z^a$, then the following estimate holds
\begin{equation}
\label{approx.est}
\norma{\delta(t)}\leq \norma{\delta(0)}e^{\rho_1 t} +
\frac{\rho_2}{\rho_1}\tond{e^{\rho_1 t} -1},\qquad\qquad t\in[-T,T]\ .
\end{equation}
\end{lemma}

\proof The difference $\delta(t)$ is solution of the differential equation
\begin{displaymath}
\dot\delta = A\delta + \tond{P(z^a(t)+\delta) - P(z^a(t))}
+ R(z^a(t))\ ,
\end{displaymath}
by Duhamel formula one has
\begin{displaymath}
\delta(t) = e^{At}\delta(0) + e^{At}\int_0^t e^{-As}
\left[\left(P(z^a(s)+\delta) - P(z^a(s))\right)+  R(z^a(s)) \right] ds
\end{displaymath}
Using Lagrange mean value theorem to estimate $P(z^a(s)+\delta) -
P(z^a(s))$ and the fact that $A$ is unitary one has
\begin{displaymath}
\norma{\delta(t)} \leq \norma{\delta(0)} + \int_0^t{\quadr{3\rho_1
\norma{\delta(s)} + \rho_2 }ds},
\end{displaymath} 
which fulfills \eqref{Gronw.hyp.1}, from which the thesis follows.\qed

\vskip 24 pt

%\bibliography{../../libro/biblio}
%\bibliographystyle{amsalpha}

\providecommand{\bysame}{\leavevmode\hbox to3em{\hrulefill}\thinspace}
\providecommand{\MR}{\relax\ifhmode\unskip\space\fi MR }
\providecommand{\MRhref}[2]{%
  \href{http://www.ams.org/mathscinet-getitem?mr=#1}{#2}
}
\providecommand{\href}[2]{#2}

\bigskip

\end{document}